\documentclass{amsart}

\usepackage{amsthm}
\usepackage{amssymb}
\usepackage{amsmath}
\usepackage{amsfonts}
\usepackage{textcomp}
\usepackage[english]{babel}
\usepackage[cp1250]{inputenc}
\usepackage{textcomp}
\usepackage{wasysym}
\usepackage{stmaryrd}
\usepackage{esint}
\usepackage[all]{xy}
\usepackage{graphicx}

\newtheorem{tw}{Theorem}[section]
\newtheorem{dfn}[tw]{Definition}
\newtheorem{uw}[tw]{Remark}
\newtheorem{prz}[tw]{Example}

\newtheorem{lem}[tw]{Lemma}
\newtheorem{stw}[tw]{Proposition}
\newtheorem{wn}[tw]{Corollary}

\newtheorem*{dd}{Proof}
\newtheorem*{sd}{Idea of proof}
\newtheorem*{ak}{Acknowledgements}

\let\olddfn\dfn
\renewcommand{\dfn}{\olddfn\normalfont}

\let\oldlem\lem
\renewcommand{\lem}{\oldlem\normalfont}

\let\oldstw\stw
\renewcommand{\stw}{\oldstw\normalfont}

\let\olduw\uw
\renewcommand{\uw}{\olduw\normalfont}

\let\oldwn\wn
\renewcommand{\wn}{\oldwn\normalfont}

\let\oldprz\prz
\renewcommand{\prz}{\oldprz\normalfont}

\let\olddd\dd
\renewcommand{\dd}{\olddd\normalfont}

\let\oldsd\sd
\renewcommand{\sd}{\oldsd\normalfont}

\let\oldak\ak
\renewcommand{\ak}{\oldak\normalfont}

\newcommand{\twopartdef}[4]
{
\left\{
		\begin{array}{ll}
			#1 & \mbox{if } #2 \\
			#3 & #4
		\end{array}
	\right.
}

\author{Wojciech G\'{o}rny}

\address{W. G\'{o}rny: Faculty of Mathematics, Informatics and Mechanics, University of Warsaw, Warsaw, Poland.}

\email{wgorny@student.uw.edu.pl}


\title{Planar least gradient problem: existence, regularity and anisotropic case}


\begin{document}

\begin{abstract}
We show existence of solutions to the least gradient problem on the plane for boundary data in $BV(\partial\Omega)$. We also provide an example of a function $f \in L^1(\partial\Omega) \backslash (C(\partial\Omega) \cup BV(\partial\Omega))$, for which the solution exists. We also show non-uniqueness of solutions even for smooth boundary data in the anisotropic case for a nonsmooth anisotropy. We additionally prove a regularity result valid also in higher dimensions.
\end{abstract}

\maketitle


\section{Introduction}

Many papers, including \cite{SWZ}, \cite{MNT}, \cite{MRL}, \cite{GRS} describe the least gradient problem, i.e. a problem of minimalization

\begin{equation*}
\min \{ \int_\Omega |Du|, \quad u \in BV(\Omega), \quad u|_{\partial\Omega} = f  \},
\end{equation*}
where we may impose certain conditions on $\Omega$, $f$ and use different approaches to the boundary condition. In \cite{SWZ} $f$ is assumed to be continuous and the boundary condition is in the sense of traces. They also impose a set of geometrical conditions on $\Omega$, which are satisfied by strictly convex sets; in fact, in dimension two they are equivalent to strict convexity. The authors of \cite{MNT} also add a positive weight. Another approach is presented in \cite{MRL}, where boundary datum belongs to $L^1(\partial\Omega)$, but the boundary condition is understood in a weaker sense.

Throughout this paper $\Omega \subset \mathbb{R}^N$ shall be an open, bounded, strictly convex set with Lipschitz (or $C^1$) boundary. The boundary datum $f$ will belong to $L^1(\partial \Omega)$ or $BV(\partial\Omega)$. We consider the following minimalization problem called the least gradient problem (for brevity denoted by LGP):

\begin{equation}\label{zagadnienie}
\inf \{ \int_\Omega |Du|, \quad u \in BV(\Omega), \quad Tu = f \},
\end{equation}
where $T$ denotes the trace operator $T: BV(\Omega) \rightarrow L^1(\partial\Omega)$. Even existence of solutions in this sense is not obvious, as the functional

\begin{equation*}
F(u) = \twopartdef{\int_\Omega |Du|}{ u \in BV(\Omega) \text{ and } Tu = f;}{+\infty}{\text{otherwise}}
\end{equation*}
is not lower semicontinuous with respect to $L^1$ convergence. In fact, in \cite{ST} the authors have given an example of a function $f$ without a solution to corresponding least gradient problem. It was a characteristic function of a certain fat Cantor set. Let us note that it does not lie in $BV(\partial\Omega)$.

There are two possible ways to deal with Problem \eqref{zagadnienie}. The first is the relaxation of the functional $F$. Such reformulation and its relationship with the original statement is considered in \cite{MRL} and \cite{Maz}. Another way is to consider when Problem \eqref{zagadnienie} has a solution in the classical sense and what is its regularity. This paper uses the latter approach.

The main result of the present paper is giving a sufficient condition for existence of solutions of the least gradient problem on the plane. It is given in the following theorem, which will be later proved as Theorem \ref{tw:istnienie}:

\begin{tw}
Let $\Omega \subset \mathbb{R}^2$ be an open, bounded, strictly convex set with $C^1$ boundary. Then for every $f \in BV(\partial\Omega)$ there exists a solution of LGP for $f$.
\end{tw}

Obviously, this condition is not necessary; the construction given in \cite{SWZ} does not require the boundary data to have finite total variation. We also provide an example of a function $f \in L^1(\Omega) \backslash (C(\partial\Omega) \cup BV(\partial\Omega))$, for which the solution exists, see Example \ref{ex:cantor}.

Another result included in this article provides a certain regularity property. Theorem \ref{tw:rozklad} asserts existence of a decomposition of a function of least gradient into a continuous and a locally constant function. It is not a property shared by all BV functions, see \cite[Example 4.1]{AFP}.

\begin{tw}
Let $\Omega \subset \mathbb{R}^N$, where $N \leq 7$, be an open, bounded, strictly convex set with Lipschitz boundary. Suppose $u \in BV(\Omega)$ is a function of least gradient. Then there exist functions $u_c, u_j \in BV(\Omega)$ such that $u = u_c + u_j$ and $(Du)_c = Du_c$ and $(Du)_j = Du_j$, i.e. one can represent $u$ as a sum of a continuous function and a piecewise constant function. They are of least gradient in $\Omega$. Moreover this decomposition is unique up to an additive constant.
\end{tw}

The final chapter takes on the subject of anisotropy. As it was proved in \cite{JMN}, for an anisotropic norm $\phi$ on $\mathbb{R}^N$ smooth with respect to the Euclidean norm there is a unique solution to the anisotropic LGP. I consider $p-$norms on the plane for $p \in [1, \infty]$ to show that for $p = 1, \infty$, i.e. where the anisotropy is not smooth, the solutions need not be unique even for smooth boundary data (see Examples \ref{ex:l1} and \ref{ex:linfty}), whereas for $1 < p < \infty$, when the anisotropy is smooth, Theorem \ref{tw:anizotropia} asserts that the only connected minimal surface with respect to the p-norm is a line segment, similarly to the isotropic solution.

\begin{tw}
Let $\Omega \subset \mathbb{R}^2$ be an open convex set. Let the anisotropy be given by the function $\phi(x,Du) = \| Du \|_p$, where $1 < p < \infty$. Let $E$ be a $\phi-$minimal set with respect to $\Omega$, i.e. $\chi_E$ is a function of $\phi-$least gradient in $\Omega$. Then every connected component of $\partial E$ is a line segment.
\end{tw}

\section{Preliminaries}

\subsection{Least gradient functions}

Now we shall briefly recall basic facts about least gradient functions. What we need most in this paper is the Miranda stability theorem and the relationship between functions of least gradient and minimal surfaces. For more information, see \cite{Giu}.

\begin{dfn}
We say that $u \in BV(\Omega)$ is a function of least gradient, if for every compactly supported $($equivalently: with trace zero$)$ $v \in BV(\Omega)$ we have

\begin{equation*}
\int_\Omega |Du| \leq \int_\Omega |D(u + v)|.
\end{equation*}
\end{dfn}

\begin{dfn}
We say that $u \in BV(\Omega)$ is a solution of the least gradient problem in the sense of traces $($solution of LGP$)$ for given $f \in L^1(\Omega)$, if $Tu = f$ and for every $v \in BV(\Omega)$ such that $Tv = 0$ we have 

\begin{equation*}
\int_\Omega |Du| \leq \int_\Omega |D(u + v)|.
\end{equation*}
\end{dfn}

To underline the difference between the two notions, we recall a stability theorem by Miranda:

\begin{tw}
$($\cite[Theorem 3]{Mir}$)$
Let $\Omega \subset \mathbb{R}^N$ be open. Suppose $\{ f_n \}$ is a sequence of least gradient functions in $\Omega$ convergent in $L^1_{loc}(\Omega)$ to $f$. Then $f$ is of least gradient in $\Omega$. \qed
\end{tw}

An identical result for solutions of least gradient problem is impossible, as the trace operator is not continuous in $L^1$ topology. We need an additional assumption regarding traces. A correct formulation would be:

\begin{tw}\label{stabilnosc}
Suppose $f, f_n \in L^1(\partial\Omega)$. Let $u_n$ be a solution of LGP for $f_n$, i.e. $Tu_n = f_n$. Let $f_n \rightarrow f$ in $L^1(\partial\Omega)$ and $u_n \rightarrow u$ in $L^1(\Omega)$. Assume that also $Tu = f$. Then $u$ is a solution of LGP for $f$.
\end{tw}
To deal with regularity of solutions of LGP, it is convenient to consider superlevel sets of $u$, i.e. sets of the form $\partial \{ u > t \}$ for $t \in \mathbb{R}$. It follows the the two subsequent results:

\begin{lem}\label{lem:jednoznacznoscnadpoziomic}
Suppose $u_1, u_2 \in L^1(\Omega)$. Then $u_1 = u_2$ a.e. iff for every $t \in \mathbb{R}$ the superlevel sets of $u_1$ and $u_2$ are equal, i.e. $\{ u_1 > t \} = \{ u_2 > t \}$ up to a set of measure zero. \qed
\end{lem}

\begin{tw}\label{twierdzeniezbgg}
$($\cite[Theorem 1]{BGG}$)$ \\
Suppose $\Omega \subset \mathbb{R}^N$ is open. Let $f$ be a function of least gradient in $\Omega$. Then the set $\partial \{ f > t \}$ is minimal in $\Omega$, i.e. $\chi_{\{ f > t \}}$ is of least gradient for every $t \in \mathbb{R}$. \qed
\end{tw}

It follows from \cite[Chapter 10]{Giu} that in low dimensions $(N \leq 7)$ the boundary $\partial E$ of a minimal set $E$ is an analytical hypersurface $($after modification of $E$ on a set of measure zero$)$. Thus, as we modify each superlevel set of $u$ by a set of measure zero, from Lemma \ref{lem:jednoznacznoscnadpoziomic} we deduce that the class of $u$ in $L^1(\Omega)$ does not change. After a change of representative we get that the boundary of each superlevel set of $u$ is a sum of analytical minimal surfaces; thus, we may from now on assume that we deal with such a representative. Also, several proofs are significantly simplified if we remember that in dimension two there is only one minimal surface: an interval.

\subsection{Sternberg-Williams-Ziemer construction}

In \cite{SWZ} the authors have shown existence and uniqueness of solutions of LGP for continuous boundary data and strictly convex $\Omega$ (or, to be more precise, the authors assume that $\partial \Omega$ has non-negative mean curvature and is not locally area-minimizing). The proof of existence is constructive and we shall briefly recall it. The main idea is reversing Theorem \ref{twierdzeniezbgg} and constructing almost all level sets of the solution. According to the Lemma \ref{lem:jednoznacznoscnadpoziomic} this uniquely determines the solution.

We fix the boundary data $g \in  C(\partial \Omega)$. By Tietze theorem it has an extension $G \in C(\mathbb{R}^n \backslash \Omega)$. We may also demand that $G \in BV(\mathbb{R}^n \backslash \overline{\Omega})$. Let $L_t =  (\mathbb{R}^n \backslash \Omega) \cap \{ G \geq t \}$. Since $G \in BV(\mathbb{R}^n \backslash \overline{\Omega})$, then for a.e. $t \in \mathbb{R}$ we have $P(L_t, \mathbb{R}^n \backslash \overline{\Omega}) < \infty$. Let $E_t$ be a set solving the following problems:

\begin{equation}\label{sternbergminimalnadlugosc}
\min \{ P(E, \mathbb{R}^n): E \backslash \overline{\Omega} = L_t \backslash \overline{\Omega} \},
\end{equation}

\begin{equation*}
\max \{ |E|: E \text{ is a minimizer of \eqref{sternbergminimalnadlugosc}} \}.
\end{equation*}

Let us note that both of these problems have solutions; let $m \geq 0$ be the infimum in the first problem. Let $E_n$ be a sequence of sets such that $P(E_n, \Omega) \rightarrow m$. By compactness of unit ball in $BV(\Omega)$ and lower semicontinuity of the total variation we obtain $\chi_{E_{n_k}} \rightarrow \chi_E$, where

\begin{equation*}
m \leq P(E, \Omega) \leq P(E_n, \Omega) \rightarrow m.
\end{equation*}
Take $M \leq |\Omega|$ be the supremum in the second problem. Take a sequence o sets $E_n$ such that $|E_n| \rightarrow M$. Then on some subsequence $\chi_{E_{n_k}} \rightarrow \chi_E$, and thus

\begin{equation*}
M \geq |E| \geq |E_n| - |E_n \triangle E| = |E_n| - \|\chi_{E_n} - \chi_E \|_1 \rightarrow M - 0.
\end{equation*}

Then we can show existence of a set $T$ of full measure such that for every $t \in T$ we have $\partial E_t \cap \partial \Omega \subset g^{-1}(t)$ and for every $t,s \in T$, $s < t$ the inclusion $E_t \subset \subset E_s$ holds. It enables us to treat $E_t$ as superlevel sets of a certain function; we define it by the following formula:

\begin{equation*}
u(x) = \sup \{t \in T: x \in \overline{E_t \cap \Omega} \}.
\end{equation*}

It turns out that $u \in C(\overline{\Omega}) \cap BV(\Omega)$ and $u$ is a solution to LGP for $g$. Moreover $| \{ u \geq t \} \triangle (\overline{E_t \cap \Omega})| = 0$ for a.e. $t$. Uniqueness proof is based on a maximum principle.

In the existence proof in chapter $4$ we are going to use a particularly simple case of the construction. Suppose $\Omega \subset \mathbb{R}^2$ and that $f \in C^1(\partial\Omega)$. Firstly, let us notice that we only have to construct the set $E_t$ for almost all $t$. Secondly, we recall that in dimension $2$ the only minimal surfaces are intervals; thus, to find the set $E_t$, let us fix $t$ and look at the preimage $g^{-1}(t)$. We connect its points with intervals with sum of their lengths as small as possible. It can cause problems, for example if we take $t$ to be a global maximum of the function; thus, let us take $t$ to be a regular value (by Sard theorem almost all values are regular), so the preimage $f^{-1}(t)$ is a manifold. In dimension $2$ this means that the preimage contains finitely many points, because $f$ is Lipschitz and $\partial\Omega$ is compact. As the derivative at every point $p \in f^{-1}(t)$ is nonzero, there is at least one interval in $\partial E_t$ ending in $p$. As is established later in Proposition \ref{slabazasadamaksimum2}, by minimality of $\partial E_t$ there can be at most one, so there is exactly one interval in $\partial E_t$ ending in every $p \in f^{-1}(t)$.

A typical example for the construction, attributed to John Brothers, is to let $\Omega = B(0,1)$ and take the boundary data to be (in polar coordinates, for fixed $r = 1$) the function $f: [0, 2\pi) \rightarrow \mathbb{R}$ given by the formula $f(\theta) = \cos(2 \theta)$; see \cite[Example 2.7]{MRL} or \cite[Example 3.6]{SZ}.

\subsection{BV on a one-dimensional compact manifold}

In the general case one may attempt to define BV spaces on compact manifolds using partition of unity; such approach is presented in \cite{AGM}. It is not necessary for us; it suffices to consider one-dimensional case. Let us consider $\Omega \subset \mathbb{R}^2$ open, bounded with $C^1$ boundary. We may define on $\partial\Omega$ the Hausdorff measure, integrable functions $($which are appoximatively continuous a.e.$)$. We recall (see \cite[Chapter 5.10]{EG}) that the one-dimensional $BV$ space on the interval $(a,b) \subset \mathbb{R}$ may be described in the following way: 

\begin{equation*}
f \in BV((a,b)) \Leftrightarrow \sum |f(x_{i})-f(x_{i-1})| \leq M < \infty
\end{equation*}
for every $a < x_0 < ... < x_n < b$, where $x_i$ are points of approximate continuity of $f$. The smallest such constant $M$ turns out to be the usual total variation of $f$.

We may extend this definition to the case where we have a one-dimensional manifold diffeomorphic to an open interval if it is properly parametrized, i.e. all tangent vectors have length one. Repeating the proof from \cite{EG} we get that this definition coincides with the divergence definition. Then we extend it to the case of a one-dimensional compact connected manifold in the following way:

\begin{dfn}
We say that $f \in BV(\partial\Omega)$, if after removing from $\partial\Omega$ a point $p$ of approximate continuity of $f$ we have
$f \in BV(\partial\Omega \backslash \{ p \})$. The norm is defined to be

\begin{equation*}
\| f \|_{BV(\partial\Omega)} = \| f \|_1 + \| f \|_{BV(\partial\Omega \backslash \{ p \})}.
\end{equation*}
\end{dfn}

This definition does not depend on the choice of $p$, as in dimension one the total variation on disjoint intervals is additive, thus for different points $p_1, p_2$ we get that 

\begin{equation*}
\| f \|_{BV(\partial\Omega \backslash \{ p_1 \})} = \| f \|_{BV((p_1,p_2))} + \| f \|_{BV((p_2,p_1))} = \| f \|_{BV(\partial\Omega \backslash \{ p_2 \})},
\end{equation*}
where $(p_1,p_2)$ is an oriented arc from $p_1$ to $p_2$. Thus all local properties of $BV(\partial\Omega)$ hold; we shall recall the most important one for our considerations:

\begin{stw}\label{stw:bvdim1}
Let $E \subset \partial\Omega$ be a set of finite perimeter, i.e. $\chi_E \in BV(\partial\Omega)$. Then, if we take its representative to be the set of points of density $1$, then $\partial E = \partial^{*} E = \{ p_1, ..., p_{2n} \}$ and $P(E, \partial\Omega) = 2n$. Here $\partial^{*} E$ denotes the reduced boundary of $E$, i.e. the set where a measure-theoretical normal vector exists; see \cite[Chapter 5.7]{EG}. \qed
\end{stw}

However, some global properties need not hold. For example, the decomposition theorem $f = f_{ac} + f_j + f_s$ does not hold; consider $\Omega = B(0,1)$, $f = \arg (z)$. The main reason is that $\pi_1(\partial\Omega) \neq 0$.

\section{Regularity of least gradient functions}

In this section we are going to prove several regularity results about functions of least gradient, valid up to dimension $7$. We start with a weak form of the maximum principle and later prove a result on decomposition of a least gradient function into a continuous and jump-type part; this decomposition holds not only at the level derivatives, but also at the level of functions. We will extensively use the blow-up theorem, see \cite[Section 5.7.2]{EG}.

\begin{tw}\label{tw:blowup}
For each $x \in \partial^{*} E$ define the set $E^r = \{ y \in \mathbb{R}^N: r(y-x) + x \in E \}$ and the hyperplane $H^{-}(x) = \{ y \in \mathbb{R}^N: \nu_E (x) \cdotp (y - x) \leq 0 \}$. Then

\begin{equation*}
\chi_{E^r} \rightarrow \chi_{H^{-}(x)}
\end{equation*}
in $L^1_{loc}(\mathbb{R}^N)$ as $r \rightarrow 0$. \qed
\end{tw}

It turns out that on the plane Theorem \ref{twierdzeniezbgg} may be improved to an analogue of the maximum principle for linear equations; geometrically speaking, the linear weak maximum principle states that every level set touches the boundary.

\begin{stw}\label{slabazasadamaksimum}
$($weak maximum principle on the plane$)$ \\
Let $\Omega \subset \mathbb{R}^2$ be an open bounded set with Lipschitz boundary and suppose $u \in BV(\Omega)$ is a function of least gradient. Then for every $t \in \mathbb{R}$ the set $\partial \{ u > t \}$ is empty or it is a sum of intervals, pairwise disjoint in $\Omega$, such that every interval connects two points of $\partial \Omega$.
\end{stw}

\begin{dd}
By the argument from \cite[Chapter 10]{Giu} for every $t \in \mathbb{R}$ the set $\partial \{ u > t \}$ is a sum of intervals and $\partial \{ u > t \} = \partial^* \{ u > t \}$. Obviously $\partial \{ u > t \}$ is closed in $\Omega$. Suppose one of those intervals ends in $x \in \Omega$. Then the normal vector at $x$ is not well defined (the statement of the Theorem \ref{tw:blowup} does not hold), so $x \notin \partial^* \{ u > t \}$. Thus $x \notin \partial \{ u > t \}$, contradiction. Similarly suppose two such intervals intersect in $x \in \Omega$. Then the measure-theoretic normal vector at $x$ has length smaller then $1$, depending on the angle between the two intervals. Thus $x \notin \partial^* \{ u > t \}$, contradiction. \qed
\end{dd}

If we additionally assume that $\Omega$ is convex, then the union is disjoint also on $\partial\Omega$:

\begin{stw}\label{slabazasadamaksimum2}
Let $\Omega \subset \mathbb{R}^2$ be an open, bounded, convex set with Lipschitz boundary and suppose $u \in BV(\Omega)$ is a function of least gradient. Then for every $t \in \mathbb{R}$ the set $\partial \{ u > t \}$ is empty or it is a sum of intervals, pairwise disjoint in $\overline{\Omega}$, such that every interval connects two points of $\partial \Omega$.
\end{stw}

\begin{dd}
Suppose that at least two intervals in $\partial E_t$ end in $x \in \partial\Omega$: $\overline{xy}$ and $\overline{xz}$. We have two possibilities: there are countably many intervals in $\partial E_t$, which end in $x$, with the other end lying in the arc $\overline{yz} \subset \partial\Omega$ which does not contain $x$; or there are finitely many. The first case is excluded by the monotonicity formula for minimal surfaces, see for example \cite[Theorem 17.6, Remark 37.9]{Sim}, as from Theorem \ref{twierdzeniezbgg} $E$ is a minimal set and only finitely many connected components of the boundary of a minimal set may intersect any compact subset of $\Omega$.

In the second case we may without loss of generality assume that $\overline{xy}$ and $\overline{xz}$ are adjacent. Consider the function $\chi_{E_t}$. In the area enclosed by the intervals $\overline{xy}, \overline{xz}$ and the arc $\overline{yz} \subset \partial\Omega$ not containing $x$ we have $\chi_{E_t} = 1$ and $\chi_{E_t} = 0$ on the two sides of the triangle (or the opposite situation, which we handle similarly). Then $\chi_{E_t}$ is not a function of least gradient: the function $\widetilde{\chi_{E_t}} = \chi_{E_t} - \chi_{\Delta xyz}$ has strictly smaller total variation due to the triangle inequality. This contradicts Theorem \ref{twierdzeniezbgg}. \qed
\end{dd}

The result above is sharp. As the following example shows, we may not relax the assumption of convexity of $\Omega$.

\begin{prz}
Denote by $\varphi$ the angular coordinate in the polar coordinates on the plane. Let $\Omega = B(0,1) \backslash (\{ \frac{\pi}{4} \leq \varphi \leq \frac{3\pi}{4} \} \cup \{ 0 \}) \subset \mathbb{R}^2$, i.e. the unit ball with one quarter removed. Take the boundary data $f \in L^1(\partial \Omega)$ to be
$$ f(x,y) = \twopartdef{1}{y \geq 0}{0}{y < 0.}$$
Then the solution to the least gradient problem is the function (defined inside $\Omega$)
$$ u(x,y) = \twopartdef{1}{y \geq 0}{0}{y < 0,}$$
in particular $\partial \{ u \geq 1 \}$ consists of two horizontal line segments whose closures intersect at the point $(0,0) \in \partial\Omega$. Note that in this example the set $\Omega$ is star-shaped, but it is not convex. \qed
\end{prz}

In higher dimensions, we are going to need a result from \cite{SWZ} concerning minimal surfaces:

\begin{stw}\label{stw:sternbergpowmin} $($\cite[Theorem 2.2]{SWZ}) \\
Suppose $E_1 \subset E_2$ and let $\partial E_1, \partial E_2$ are area-minimizing in a open set $U$. Further, suppose $x \in \partial E_1 \cap \partial E_2 \cap U$. Then $\partial E_1$ and $\partial E_2$ agree in some neighbourhood of $x$. \qed
\end{stw}

\begin{tw}$($weak maximum principle$)$ \\
Let $\Omega \subset \mathbb{R}^N$, where $N \leq 7$ and suppose $u \in BV(\Omega)$ is a function of least gradient. Then for every $t \in \mathbb{R}$ the set $\partial \{ u > t \}$ is empty or it is a sum of minimal surfaces $S_{t,i}$, pairwise disjoint in $\Omega$, which satisfy $\partial S_{t,i} \subset \partial \Omega$.
\end{tw}

\begin{dd}
Let us notice, that with only subtle changes the previous proof works also in the case $N \leq 7$, i.e. when boundaries of superlevel sets are minimal surfaces.

From \cite[Chapter 10]{Giu} it follows that for $t \in \mathbb{R}$ the set $\partial \{ u > t \}$ is a sum of minimal surfaces $S_{t,i}$ and $\partial \{ u > t \} = \partial^* \{ u > t \}$. Obviously $\partial \{ u > t \}$ $($boundary in topology of $\Omega)$ is closed in $\Omega$, so $\partial S_{t,i} \cap \Omega = \emptyset$ $($boundary in topology of $\partial \{ u > t \})$; suppose otherwise. Let $x \in \partial S_{t,i} \cap \Omega$. Then in $x$ the blow-up theorem does not hold, so $x \notin \partial \{ u > t \}$, contradiction.

Now suppose that $S_{t,i}$ and $S_{t,j}$ are not disjoint in $\Omega$. Then from the Proposition \ref{stw:sternbergpowmin} applied to $E_1 = E_2 = \{ u > t \}$ we get $S_{t,i} = S_{t,j}$. \qed
\end{dd}

\begin{stw}
Let $E_1 \subset E_2$ and suppose that $E_1$ and $E_2$ are sets of locally bounded perimeter and let $x \in \partial^{*} E_1 \cap \partial^{*} E_2$. Then $\nu_{E_1}(x) = \nu_{E_2}(x)$.
\end{stw}

\begin{dd}
We are going to use the blow-up theorem (Theorem \ref{tw:blowup}). First notice that the inclusion $E_1 \subset E_2$ implies

\begin{equation*}
E_1^r = \{ y \in \mathbb{R}^N: r(y-x) + x \in E_1 \} \subset \{ y \in \mathbb{R}^N: r(y-x) + x \in E_2 \} = E_2^r.
\end{equation*}
We keep the same notation as in Theorem \ref{tw:blowup} and use it to obtain

\begin{equation*}
\chi_{H_1^-(x)} \leftarrow \chi_{E_1^r} \leq \chi_{E_1^r} \rightarrow \chi_{H_2^-(x)},
\end{equation*}
where the convergence holds in $L^1_{loc}$ topology. Thus $H_1^-(x) = H_2^-(x)$, so $\nu_{E_2}(x) = \nu_{E_2}(x)$. \qed
\end{dd}

\begin{stw}\label{stw:zbiorskokow}
For $u \in BV(\Omega)$ the structure of its jump set is as follows:

\begin{equation*}
J_u = \bigcup_{s,t \in \mathbb{Q}; s \neq t} (\partial^{*} \{ u > s \} \cap \partial^{*} \{ u > t \}).
\end{equation*}
\end{stw}

\begin{dd}
Let $x \in J_u$. By definition of $J_u$ the normal vector at $x$ is well defined. The same applies to the trace values from both sides: let us denote them by $u^{-} (x) < u^{+} (x)$. But then there exist $s,t \in \mathbb{Q}$ such that $u^{-} (x) < s < t < u^{+} (x)$, so $x \in \partial^{*} \{ u > s \} \cap \partial^{*} \{ u > t \}$.

On the other hand, let $x \in \partial^{*} \{ u > s \} \cap \partial^{*} \{ u > t \}$. From the previous proposition the normal vectors coincide, so the normal at $x$ does not depend on $t$ and we may define traces from both sides as 

\begin{equation*}
u^{+}(x) = \sup \{ t: x \in \partial^{*} \{ u > s \} \cap \partial^{*} \{ u > t \} \};
\end{equation*}
\begin{equation*}
u^{-}(x) = \inf \{ t: x \in \partial^{*} \{ u > s \} \cap \partial^{*} \{ u > t \} \}.
\end{equation*}
More precisely, the trace is uniquely determined up to a measure zero set by the mean integral property from \cite[Theorem 5.3.2]{EG}. But it holds for all $x \in \partial^{*} \{ u > s \} \cap \partial^{*} \{ u > t \}$; from the weak maximum principle this set divides $\Omega$ into two disjoint parts, $\Omega^+$ and $\Omega^-$. Let $\Omega^+$ be the part with greater values of $u$ in the neighbourhood of the cut. If $u^+(x) < \sup \{ t: x \in \partial^{*} \{ u > s \} \cap \partial^{*} \{ u > t \} \} = s$, then for sufficiently small neighbourhoods of $x$ we would have $u \geq s$, so $\fint_{B(x,r) \cap \Omega^+} |u^+(x) - u(y)| \geq |u^+(x) - s| > 0$, contradiction. The other cases are analogous. \qed
\end{dd}

\begin{stw}
Suppose $u \in BV(\Omega)$ is a least gradient function. Then $J_u = \bigcup_{k=1}^{\infty} S_k$, where $S_k$ are pairwise disjoint minimal surfaces. In addition, the trace of $u$ from both sides is constant along $S_k$; in particular the jump of $u$ is constant along $S_k$.
\end{stw}

\begin{dd}
We follow the characterisation of $J_u$ from the Proposition \ref{stw:zbiorskokow}. For every $t$ the set $\partial^{*} \{ u > t \}$ is a minimal surface. Proposition \ref{stw:sternbergpowmin} ensures that if $\partial^{*} \{ u > s \} \cap \partial^{*} \{ u > t \} \neq \emptyset$, then their intersecting connected components $S_{s,i}$, $S_{t,j}$ coincide. In particular, the trace from both sides defined as above is constant along $S_{t,j}$. Thus connected components of $J_u$ coincide with connected components of $\partial^{*} \{ u > t \}$ for some $t$, so by weak maximum principle they are minimal surfaces non-intersecting in $\Omega$ with boundary in $\partial \Omega$. As the area of each such surface is positive, there is at most countably many of them. \qed
\end{dd}

\begin{tw}\label{tw:rozklad}
Let $\Omega \subset \mathbb{R}^N$, where $N \leq 7$, be an open, bounded, strictly convex set with Lipschitz boundary. Suppose $u \in BV(\Omega)$ is a function of least gradient. Then there exist functions $u_c, u_j \in BV(\Omega)$ such that $u = u_c + u_j$ and $(Du)_c = Du_c$ and $(Du)_j = Du_j$, i.e. one can represent $u$ as a sum of a continuous function and a piecewise constant function. They are of least gradient in $\Omega$. Moreover this decomposition is unique up to an additive constant.
\end{tw}

\begin{dd}
1. From the previous theorem $J_u = \bigcup_{k=1}^{\infty} S_k$, where $S_k$ are pairwise disjoint minimal surfaces with boundary in $\partial\Omega$. The jump along each of them has a constant value $a_k$. They divide $\Omega$ into open, pairwise disjoint sets $U_i$.

\begin{figure}[h]
\includegraphics[scale = 0.15]{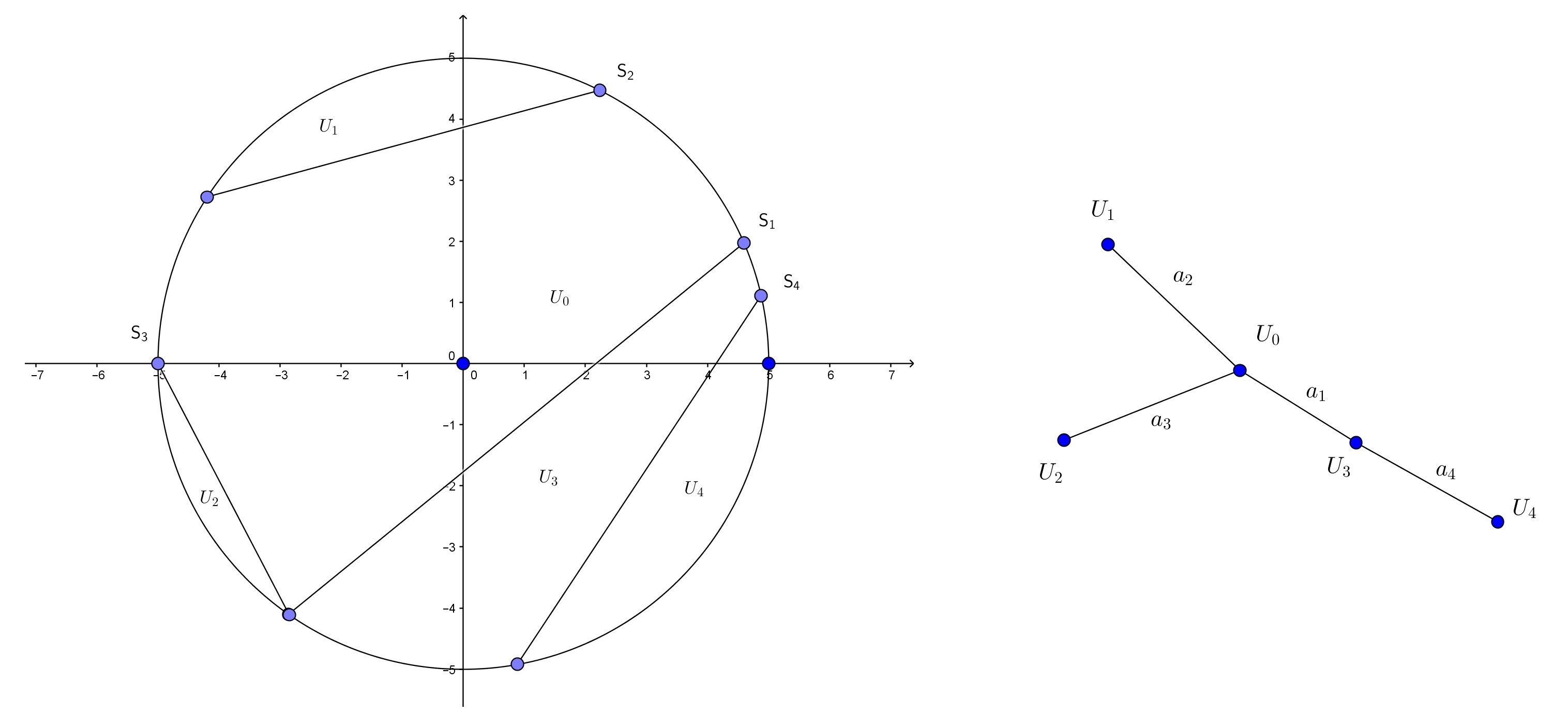}
\end{figure}

2. We define $u_j$ in the following way: let us call any of the obtained sets $U_0$. Let us draw a graph such that the sets $U_i$ are its vertices. $U_i$ and $U_j$ are connected by an edge iff $\partial U_i \cap U_j = S_k$, i.e. when they have a common part of their boundaries. To such an edge we ascribe a weight $a_k$. Example of such construction is presented on the picture above. As $S_k$ are disjoint in $\Omega$ and do not touch the boundary, then such a graph is a tree, i.e. it is connected and there is exactly one path connecting two given vertices. Thus, we define $u_j$ by the formula

\begin{equation*}
u_j(x) = \sum_{\text{path connecting } U_0 \text{ with } U_i} a_k, \text{ when } x \in U_j.
\end{equation*}
Such a function is well defined, as our graph has no cycles. It also does not depend on the choice of $U_0$ up to an additive constant (if we chose some $U_1$ instead, the function would change by a summand $\sum_{\text{path connecting } U_0 \text{ with } U_1} a_k)$. We see that $u_j \in L^1(\Omega)$ and that it is piecewise constant. \\

3. We notice that $Du_j = (Du)_j$, as $u_j$ is constant on each $U_i$, $J_{u_j} = J_u$ and the jumps along connected components of $J_u$ have the same magnitude. Thus we define $u_c = u - u_j$. We see that $(Du_c)_j = 0$. \\

4. The $u_c$, $u_j$ defined above are functions of least gradient.

Suppose that $u_j$ is not a a function of least gradient, i.e. there exists $v \in BV(\Omega)$ such that $\int_\Omega |Dv| < \int_\Omega |Du_j|$ and $Tu_j = Tv$. Then we would get

\begin{equation*}
\int_\Omega |Du| \leq \int_\Omega |D(u_c + v)| \leq \int_\Omega |Du_c| + \int_\Omega |Dv| < \int_\Omega |Du_c| + \int_\Omega |Du_j| = \int_\Omega |Du|,
\end{equation*}
where the first inequality follows from $u$ being a function of least gradient, and the last equality from measures $Du_c$ and $Du_j$ being mutually singular. The proof for $u_c$ is analogous.

5. The function $u_c$ is continuous. As $u_c$ is of least gradient, then if it isn't continuous at $x \in \Omega$, then a certain set of the form $\partial \{ u_c > t \}$ passes through $x$; otherwise $u_c$ would be constant in the neighbourhood of $x$. But in that case $u_c$ has a jump along the whole connected component of $\partial \{ u_c > t \}$ containing $x$, which is impossible as $(Du_c)_j = 0$.

6. What is left is to prove uniqueness of such a decomposition. Let $u = u_c^1 + u_j^1 = u_c^2 + u_j^2$. Changing the order of summands we obtain

\begin{equation*}
u_c^1 - u_c^2 = u_j^2 - u_j^1,
\end{equation*}
but the distributional derivative of the left hand side is a continuous measure, and the distributional derivative of the right hand side is supported on the set of zero measure with respect to $\mathcal{H}^{n-1}$, so both of them are zero measures. But the condition $Dv = 0$ implies $v = \text{const}$, so the functions $u_c^1$, $u_c^2$ differ by an additive constant. \qed
\end{dd}

\begin{prz}
In this decomposition $u_c$ isn't necessarily continuous up to the boundary. Let us use the complex numbers notation for the plane. We take $\Omega = B(1,1)$. Let the boundary values be given by the formula $f(z) = \arg(z)$. Then $u = u_c = \arg(z) \in BV(\Omega) \cap C^{\infty}(\Omega)$, but $u$ isn't continuous at $0 \in \partial\Omega$. \qed
\end{prz}

\section{Existence of solutions on the plane}

We shall prove existence of solutions on the plane for boundary data in $BV(\partial\Omega)$. We are going to use approximations of the solution in strict topology. Proposition \ref{podciagzbiezny} will ensure us that existence of convergent sequences of approximations in $L^1$ topology is not a problem; Theorem \ref{tw:scislazb} will upgrade it to strict convergence. The Miranda stability theorem (Theorem \ref{stabilnosc}) ends the proof. Later, we shall see an example of a discontinuous function $f$ of infinite total variation such that the solution to the LGP exists.

\begin{stw}\label{podciagzbiezny}
Suppose $f_n \rightarrow f$ in $L^1(\partial\Omega)$. $u_n$ are solutions of LGP for $f_n$. Then $u_n$ has a convergent subsequence, i.e. $u_{n_k} \rightarrow u$ in $L^1(\Omega)$.
\end{stw}

\begin{dd}
As the trace operator is a surjection, by the Open Mapping Theorem it is open. Let us fix $\widetilde{f} \in BV(\Omega)$ such that $T\widetilde{f} = f$ and a sequence of positive numbers $\varepsilon_n \rightarrow 0$. Then by continuity and openness of $T$ the image of a ball $B(\widetilde{f}, \varepsilon_n)$  contains a smaller ball $B(T\widetilde{f}, \delta_n)$ for another sequence of positive numbers $\delta_n \rightarrow 0$. As $f_n \rightarrow f$ in $L^1(\partial\Omega)$, there exists a subsequence $f_{n_k}$ such that $f_{n_k} \in B(f, \delta_n) = B(T\widetilde{f}, \delta_n)$, so the set $T^{-1}(f_{n_k})$ is non-empty; there exists a preimage of $f_{n_k}$ by $T$ in $B(\widetilde{f}, \varepsilon_n)$. Let us call it $\widetilde{f_n}$. Obviously $\widetilde{f_n} \rightarrow \widetilde{f}$ in $BV(\Omega)$.

Thus, after possibly passing to a subsequence, there exist functions $\widetilde{f_n}, \widetilde{f}$ such that $\widetilde{f_n} \rightarrow \widetilde{f}$ in $BV(\Omega)$ and $T\widetilde{f_n} = f_n, T\widetilde{f} = f$.
Now we may proceed as in \cite[Proposition 3.3]{HKLS}. Let us estimate from above the norm of $\|u_n - \widetilde{f_n}\|_{BV}$:

\begin{multline*}
\|u_n - \widetilde{f_n}\|_{BV} = \|u_n - \widetilde{f_n}\|_1 + \int_\Omega |D(u_n - \widetilde{f_n})| \leq (C + 1) \int_\Omega |D(u_n - \widetilde{f_n})| \leq \\
\leq (C + 1) (\int_\Omega |Du_n| + \int_\Omega |D\widetilde{f_n}|) \leq 2(C + 1) \int_\Omega |D\widetilde{f_n}| \leq M < \infty
\end{multline*}
where the inequalities follow from Poincaré inequality $($as $u_n - f_n$ has trace zero$)$, triangle inequality and the fact that $u_n$ is solution of LGP for $f_n$. The common bound follows from convergence of $\widetilde{f_n}$. 

Thus, by compactness of the unit ball of $BV(\Omega)$ in $L^1(\Omega)$ we get a convergent subsequence $u_{n_k} - \widetilde{f_{n_k}} \rightarrow v$ in $L^1(\Omega)$. But $\widetilde{f_n} \rightarrow \widetilde{f}$ in $BV(\Omega)$, so as well in $L^1(\Omega)$; thus $u_{n_k} \rightarrow v + \widetilde{f} = u$ in $L^1(\Omega)$. \qed
\end{dd}

We are going to need three lemmas. The first two are straightforward and their proofs can be found as a step in the proof of co-area formula, see \cite[Section 5.5]{EG}. The third one is a convenient version of Fatou lemma.

\begin{lem}\label{lem:zb1}
Let $f_n \rightarrow f$ in $L^1(\Omega)$. Then there exists a subsequence $f_{n_k}$ such that $\chi_{\{ f_{n_k} \geq t \}} \rightarrow \chi_{\{ f \geq t \}}$ in $L^1(\Omega)$ for a.e. $t$. \qed
\end{lem}

\begin{lem}\label{lem:zb2}
Suppose $\chi_{\{ f_{n} \geq t \}} \rightarrow \chi_{\{ f \geq t \}}$ in $L^1(\Omega)$ for a.e. $t$. Then $f_n \rightarrow f$ in $L^1_{loc}(\Omega)$. If additionally $f, f_n$ form a bounded family in $L^{\infty}(\Omega)$, then this covergence holds also in $L^1(\Omega)$. \qed
\end{lem}

\begin{lem}\label{lem:zb3}
Suppose that $g, g_n \geq 0$. If additionally $g \leq \liminf g_n$ a.e. and $\lim \int_\Omega g_n \, dx = \int_\Omega g \, dx < \infty$, then $g_n \rightarrow g$ in $L^1(\Omega)$.
\end{lem}

\begin{dd}
Let $f_+ = \max(f, 0)$ and $f_- = \max(-f, 0)$. Let us note that
\begin{equation*}
\int_\Omega |g - g_n| = \int_\Omega (g - g_n)_+ + \int_\Omega (g - g_n)_-
\end{equation*}
and 

\begin{equation*}
0 \leftarrow \int_\Omega (g - g_n) = \int_\Omega (g - g_n)_+ - \int_\Omega (g - g_n)_-,
\end{equation*}
so it suffices to prove that $\int_\Omega (g - g_n)_+ \rightarrow 0$ to show that $g_n \rightarrow g$ in $L^1(\Omega)$. Now let us see what happens to (well defined) upper limit of the sequence $\int_\Omega (g - g_n)_+$:

\begin{equation*}
0 \leq \limsup \int_\Omega (g - g_n)_+ \leq \int_\Omega \limsup (g - g_n)_+ = \int_\Omega \limsup \max(g - g_n, 0) =
\end{equation*}
\begin{equation*}
= \int_\Omega \max(g + \limsup (- g_n), 0) = \int_\Omega \max(g - \liminf g_n), 0) = \int_\Omega 0 = 0.
\end{equation*}
where inequality follows from the (inverse) Fatou lemma: by definition $0 \leq (g - g_n)_+ \leq g$, and $g$ is integrable, so we can apply the Fatou lemma. To prove equalities we use the fact that $\limsup (- g_n) = - \liminf g_n$ and the assumption that $g \leq \liminf g_n$ a.e. Thus $\int_\Omega (g - g_n)_+ \rightarrow 0$, so $g_n \rightarrow g$ in $L^1(\Omega)$.
\end{dd}

\begin{tw}\label{tw:scislazb}
Let $\Omega \subset \mathbb{R}^2$ be an open, bounded, strictly convex set with $C^1$ boundary and suppose $f \in BV(\partial\Omega)$. Let $f_n \rightarrow f$ strictly in $BV(\partial\Omega)$, where $f_n$ are smooth. Denote the unique solution of LGP for $f_n$ by $u_n$. Then on some subsequence $u_{n_k}$ we have strict convergence in $BV(\Omega)$ to a function $u \in BV(\Omega)$. In particular $Tu = f$.
\end{tw}

\begin{dd}
1. As we have $f_n \rightarrow f$ strictly in $BV(\partial\Omega)$, we by definition also convergence in $L^1(\partial\Omega)$. Thus, by Lemma \ref{lem:zb1}, after possibly passing to a subsequence we have convergence $\chi_{\{ f_n \geq t \}} \rightarrow \chi_{\{ f \geq t \}}$ for a.e. $t$.

2. By co-area formula
\begin{equation*}
\int_{\partial\Omega} |Df_n| = \int_{\mathbb{R}} P(E_t^n, \partial\Omega) \, dt \rightarrow \int_{\mathbb{R}} P(E_t, \partial\Omega) dt = \int_{\partial \Omega} |Df|,
\end{equation*}
and lower semicontinuity of the total variation gives us $P(E_t, \partial\Omega) \leq \liminf P(E_t^n, \partial\Omega) < \infty$ for a.e. $t$. We observe that the conditions in Lemma \ref{lem:zb3} are fulfilled and we obtain convergence $P(E_t^n, \partial\Omega) \rightarrow P(E_t, \partial\Omega)$ in $L^1(\mathbb{R})$, so after possibly passing to a subsequence we have pointwise convergence for a.e. $t$. Consequently $\chi_{\{ f_n \geq t \}} \rightarrow \chi_{\{ f \geq t \}}$ strictly in $BV(\partial\Omega)$.

3. As $\partial\Omega \in C^1$ and $f_n \in C^1(\partial\Omega)$, then by Sard theorem the set $\mathcal{T}$ of such $t$, which are regular values for all $f_n$, is of full measure. Recalling the Sternberg-Williams-Ziemer construction we get that for every $t \in \mathcal{T}$ every point of $\partial E_t^n \cap \partial\Omega$ is an end of at least one interval; according to Proposition \ref{slabazasadamaksimum2} it is an end of exactly one interval.

4. From now on it is necessary that we are in dimension $N = 2$. Let $t \in \mathcal{T}$. As $\partial\Omega$ is one-dimensional, then $P(E_t^n, \partial\Omega) \in \mathbb{N}$ and $D \chi_{\{ f_n \geq t \}}$ is a sum $\sum_{i = 1}^M \pm \delta_{x_i}$. Furthermore, by Proposition \ref{stw:bvdim1} there exists a representative of the set $E_t^n$, which is a sum of closed arcs between consecutive points $x_i$. By Lemma \ref{lem:jednoznacznoscnadpoziomic} we can change all representatives of the sets $E_t^n$ not changing $f_n$ itself. We do the same for $E_t$. As $f_n$ are smooth functions, such form of $E_t^n$ follows directly from their smoothness; this needs not be the case for $E_t$.

5. As $\chi_{\{ f_n \geq t \}} \rightarrow \chi_{\{ f \geq t \}}$ strictly, then for sufficiently large $n$ $P(E_t^n, \partial\Omega) = P(E_t, \partial\Omega)$. What is more, their derivatives converge in weak* topology; but we have an exact representation of those derivative. This gives us convergence $x_i^n \rightarrow x_i$ for every $i$.

6. We apply the Sternberg-Williams-Ziemer construction to the sequence $f_n$. The set $\partial E_t^n$ is a sum of intervals, disjoint in $\Omega$, connecting certain pairs of points among $x_i^n$. By definition of $\mathcal{T}$ every point of $\partial E_t^n \cap \partial\Omega$ is an end of exactly one interval. This gives us convergence $\chi_{\{ u_n \geq t \}} \rightarrow \chi_{\{ u \geq t \}}$ w $L^1(\Omega)$ for a.e. $t$. Because of continuity of the metric in $\mathbb{R}^2$ we get $P(E_t^n, \Omega) = \sum \| x_i^n - x_j^n \| \rightarrow \sum \| x_i - x_j \| = P(E_t, \Omega)$.

7. Let us see that $P(E_t^n, \Omega) \leq P(\Omega, \mathbb{R}^N)$. Indeed, $\partial E_t^n$ is a sum of intervals, disjoint in $\Omega$, connecting certain pairs of points among $x_i^n$. If we choose a different connection between them, for example by drawing a full convex polygon with vertices in $x_i^n$, by minimality of $\partial E_t^n$ the polygon has a larger perimeter. If we use arcs on $\partial\Omega$ instead, the perimeter would be even larger, as intervals are minimal surfaces in $\mathbb{R}^2$.

8. Since the functions $\chi_{\{ u_n \geq t \}}$ converge in $L^1(\Omega)$ for a.e. $t$ to $\chi_{\{ u \geq t \}}$, then by Lemma \ref{lem:zb2} we have convergence $u_n \rightarrow u$ in $L^1(\Omega)$. Furthermore in step $6$ we proved convergence $P(E_t^n, \Omega) \rightarrow P(E_t, \Omega)$ for a.e. $t$, so by dominated convergence theorem (by step $7$ this sequence is bounded) we have convergence $P(E_t^n, \Omega) \rightarrow P(E_t, \Omega)$ in $L^1(\mathbb{R})$. By co-area formula $\int_\Omega |Du_n| \rightarrow \int_\Omega |Du|$, which gives that $u_n \rightarrow u$ strictly in $BV(\Omega)$.
\end{dd}

\begin{tw}\label{tw:istnienie}
Let $\Omega \subset \mathbb{R}^2$ be an open, bounded, strictly convex set with $C^1$ boundary. Then for every $f \in BV(\partial\Omega)$ there exists a solution of LGP for $f$.
\end{tw}

\begin{dd}
For each $f \in BV(\partial\Omega)$ we can find a sequence $f_n$ of class $C^{\infty}(\partial\Omega)$ strictly convergent to $f$. Let $u_n$ be solutions of LGP for $f_n$. Then after possibly passing to subsequence we have that $u_n \rightarrow u$ strictly in $BV(\Omega)$; but the trace is continuous in the strict topology, so $Tu = f$. Thus by Miranda stability theorem (Theorem \ref{stabilnosc}) we get that $u$ is a solution of LGP for $f$.
\end{dd}

\begin{prz}
Take $\Omega = B(0,1)$. As we know from \cite{ST}, when $f$ is a characteristic function of a certain fat Cantor set, then the least gradient problem has no solution. Thus, we would expect that if we approximated the boundary function and constructed solutions of LGP for the approximation, then the trace of the limit would be incorrect. To settle this, let $f_n$ be a function of the $n-$th stage of the Cantor set construction. Then $u_n \rightarrow 0$ in $L^1(\Omega)$:

Let $f_0(\theta) = \chi_{[0,1]}$. We construct $f_1$ by removing from the middle of $[0,1]$ an interval of length $2^{-2}$, i.e. $f_1 = \chi_{[0,3/8] \cup [5/8,1]}$. In the second stage we remove from the middle of both intervals an interval of length $2^{-4}$ and obtain $f_2 = \chi_{[0,5/32] \cup [7/32, 3/8] \cup [5/8,25/32] \cup [27/32, 1]}$. During the $n-$th stage of construction we remove an interval of length $2^{-2n}$ from the middle of all existing $2^{n-1}$ intervals.

Let us see what is the length of all such intervals. Let $a_n$ be the length of an interval at the $n-$th stage of construction. Then $a_{n} = \frac{a_{n-1}}{2} - \frac{1}{2^{2n+1}}$. As $a_0 = 1$, we obtain a direct formula $a_n = \frac{2^n + 1}{2^{2n+1}}$.

Now we take the fat Cantor set to be on the circle, i.e. the interval [0,1] corresponds to angles measured in radians. On the rest of the circle we set the function $f$ to be $0$.

\begin{figure}[h]
\includegraphics[scale = 0.15]{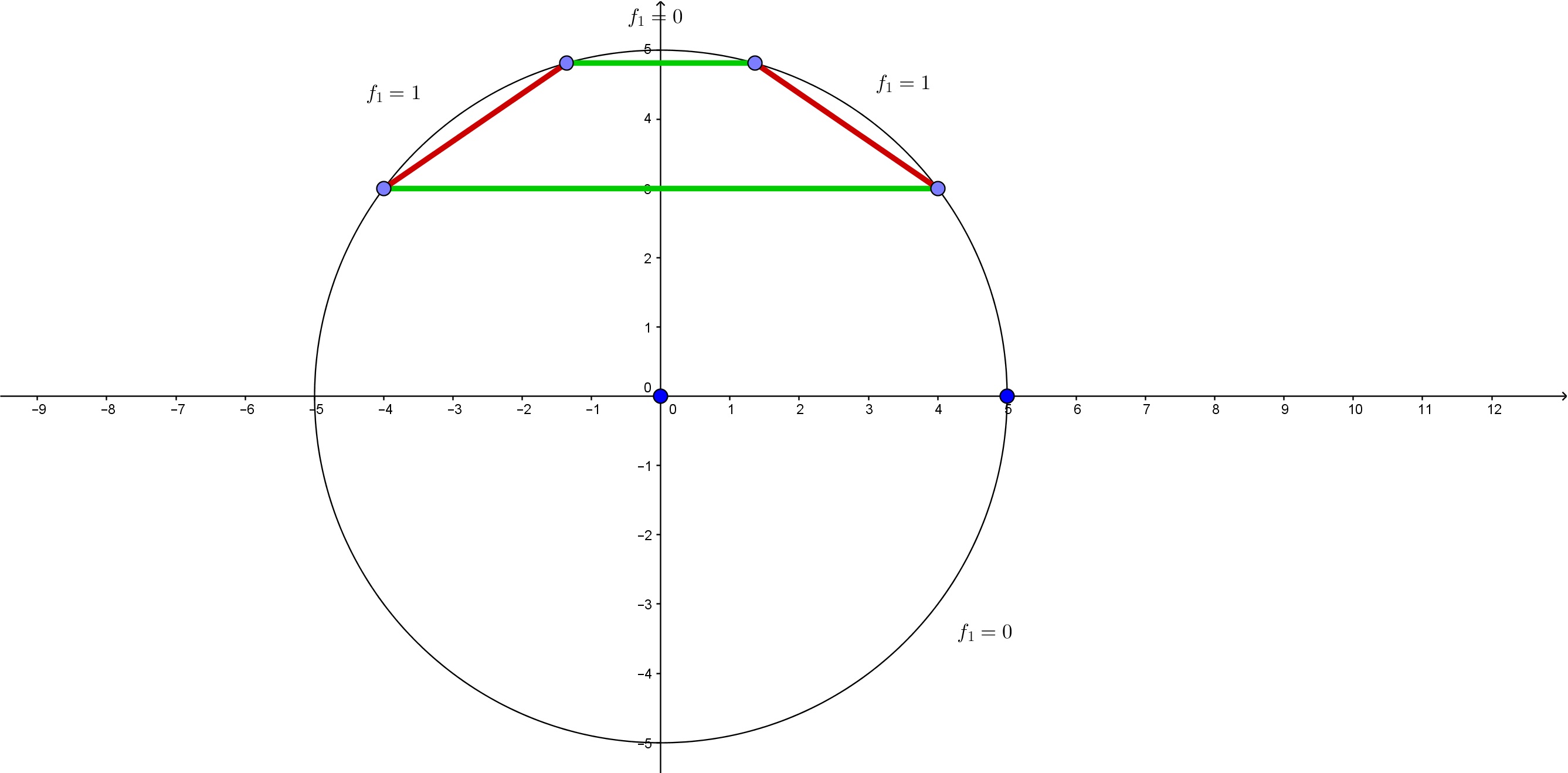}
\end{figure}

Let us compare at every stage of construction the sum of lengths of the red intervals and the green ones. After trygonometric considerations we have to check the following inequality:

\begin{equation}\label{ineq:cantor}
\sqrt{1 - \cos(a_n)} + \sqrt{1 - \cos(\frac{1}{2^{2n}})} > 2 \sqrt{1 - \cos(a_{n+1})}.
\end{equation}
Substitute $x = 2^{-n}$ and use the direct formula for $a_n$. It changes to the inequality

\begin{equation*}
g(x) = \sqrt{1 - \cos(\frac{x(x+1)}{2})} + \sqrt{1 - \cos(x^2)} - 2 \sqrt{1 - \cos(\frac{x(x+2)}{8})} > 0.
\end{equation*}
But $g$ satisfies $g(0) = 0$ and its derivative is positive on $(0,1)$, so $g > 0$ on $(0,1)$, thus the inequality holds for all $n$. Thus, as on every stage of construction the sides of the trapezoid are shorter than the bases. It means that the solution of LGP for $f_n$ takes value $0$ on the trapezoid (as we minimize $P(E_t, B(0,1))$ for $t \in (0,1)$). In the next stage of construction the value on this trapezoid will remain zero and we will make the same reasoning on two adjacent smaller trapezoids. From the construction of Cantor set the sequence $u_n$ is nonincreasing and for every point $x$ inside the circle at a sufficiently large stage of construction we would have $u_n(x) = 0$. Thus $u_n \rightarrow 0$ a.e.; but it is bounded from above by $1$, so the convergence holds also in $L^1(\Omega)$.
\end{prz}

\begin{prz}\label{ex:cantor}
Let us make a slight change to the previous example: consider another fat Cantor set. More precisely, take a set almost of full measure such that the inequality \eqref{ineq:cantor} holds in the opposite direction; it is possible due to the triangle inequality. Thus at every stage of construction it is more efficient (minimizing lengths of level sets) to remove $2^{n-1}$ curvilinear triangles from the set $\{ u_n = 1 \}$ than to repeat the above construction, i.e. add trapezoids to the set $\{ u_n = 0 \}$. Thus at every stage of construction the set $\{ u_n = 1 \}$ will be a sum of trapezoids mentioned before, so the trace of $u$ equals $f$. Also $u_n \rightarrow u$ in $L^1(\Omega)$, as it converges a.e. Thus we obtained that there exists a solution to LGP for a certain discontinuous $f \notin BV(\partial\Omega)$.
\end{prz}

\section{Anisotropic case}

This section is devoted to the anisotropic least gradient problem. We discuss $l^p$ norms on the plane for $p \in [1, \infty]$. We prove a non-uniqueness result for $p = 1, \infty$ and discuss how the solutions look like for $p \in (1, \infty)$. We shall use the notation introduced in \cite{Maz}.

\begin{dfn}
A continuous function  $\phi: \overline{\Omega} \times \mathbb{R}^n \rightarrow [0, \infty)$ is called a metric integrand, if it satisfies the following conditions: \\
\\
$(1)$ $\phi$ is convex with respect to the second variable for a.e. $x \in \overline{\Omega}$; \\
$(2)$ $\phi$ is homogeneous with respect to the second variable, i.e.

\begin{equation*}
\forall x \in \overline{\Omega}, \quad \forall \xi \in \mathbb{R}^n, \quad \forall t \in \mathbb{R} \quad \phi(x, t \xi) = |t| \phi(x, \xi);
\end{equation*}
$(3)$ $\phi$ is bounded in $\Omega$, i.e.

\begin{equation*}
\exists \Gamma > 0 \quad \forall x \in \overline{\Omega}, \quad \forall \xi \in \mathbb{R}^n \quad 0 \leq \phi(x, \xi) \leq \Gamma |\xi|.
\end{equation*}
\\
$(4)$ $\phi$ is elliptic in $\Omega$, i.e.

\begin{equation*}
\exists \lambda > 0 \quad \forall x \in \overline{\Omega}, \quad \forall \xi \in \mathbb{R}^n \quad \lambda |\xi| \leq \phi(x, \xi).
\end{equation*}
\end{dfn}

\begin{uw}
These conditions are sufficient for most of the cases considered in scientific work: they are satisfied for the classical LGP, i.e. $(\phi(x, \xi) = |\xi|)$, as well as for the $l_p$ norms, $p \in [1, \infty]$ and for weighted LGP considered in \cite{JMN}: a function $\phi(x, \xi) = g(x) |\xi|$, where $g \geq c > 0$.
\end{uw}

\begin{dfn}
The polar function of $\phi$ is $\phi^0: \overline{\Omega} \times \mathbb{R}^N \rightarrow [0, \infty)$ defined as

\begin{equation*}
\phi^0 (x, \xi^*) = \sup \, \{ \langle \xi^*, \xi \rangle : \xi \in \mathbb{R}^N, \phi(x, \xi) \leq 1 \}.
\end{equation*} 
\end{dfn}

\begin{dfn}
Let 
\begin{equation*}
\mathcal{K}_\phi(\Omega) = \{ \mathbf{z} \in X_\infty(\Omega) : \phi^0(x,\mathbf{z}(x)) \leq 1 \text{ for a.e. } x \in \Omega, \quad [\mathbf{z}, \nu] = 0 \}.
\end{equation*}
For a given function $u \in L^1(\Omega)$ we define its $\phi-$total variation in $\Omega$ by the formula (another notation used in the literature is $\int_\Omega \phi(x, Du)$):

\begin{equation*}
\int_\Omega |Du|_\phi = \sup \, \{ \int_\Omega u \, \mathrm{div} \, \mathbf{z} \, dx : \mathbf{z} \in \mathcal{K}_\phi(\Omega) \}.
\end{equation*}
If $\int_\Omega |Du|_\phi < \infty$, we say that $u \in BV_\phi(\Omega)$. If $\phi$ is a metric integrand, by properties $(3)$ and $(4)$ we have that $\lambda \int_\Omega |Du| \leq \int_\Omega |Du|_\phi \leq \Gamma \int_\Omega |Du|$, so $BV_\phi(\Omega) = BV(\Omega)$. We also know (\cite[Chapter 3]{AB}) that when $\phi$ is continuous and elliptic in $\Omega$, then in the definition of $\mathcal{K}_\phi(\Omega)$ we can replace the condition $[\mathbf{z}, \nu] = 0$ with a demand that $\mathbf{z} \in C_c^1(\Omega)$, so we recover the classical definition.
\end{dfn}

\begin{uw}\label{lem:scisaniz}
When $\phi$ is continuous and elliptic in $\Omega$, then similarly to the classical case (\cite[Chapter 4]{AB}) we recover lower semicontinuity of the $\phi-$total variation, the notion of $\phi-$perimeter of a set and the co-area formula. We also recover the approximation by $C^\infty$ functions in the strict topology, even in the strong form proved by Giusti in \cite[Corollaries 1.17, 2.10]{Giu}: let $v \in BV_\phi(\Omega)$, $Tv = f$. Then there exists a sequence of $C^\infty$ functions $v_n$ such that $v_n \rightarrow v$ strictly in $BV_\phi(\Omega)$ such that $Tv = f$. \qed
\end{uw}

For an explicit use we shall need the following integral representation (\cite{AB}, \cite{JMN}):

\begin{stw}\label{stw:repcalkowa}
Let $\varphi: \overline{\Omega} \times \mathbb{R}^N \rightarrow \mathbb{R}$ be a metric integrand. Then we have an integral representation:

\begin{equation*}
\int_\Omega |Du|_\phi = \int_\Omega \phi(x, \nu^u(x)) \, |Du|,
\end{equation*}
where $\nu^u$ is the Radon-Nikodym derivative $\nu^u = \frac{d Du}{d |Du|}$. In particular, if $E \subset \Omega$ and $\partial E$ is sufficiently smooth (at least $C^1$), then we have a representation

\begin{equation*}
P_\phi(E, \Omega) = \int_\Omega \phi(x, \nu_E) \, d \mathcal{H}^{n-1},
\end{equation*}
where $\nu_E$ is the external normal to $E$. \qed
\end{stw}

\begin{dfn}
For $p \in [1, \infty)$ we define the $p-$th norm of a vector on the plane by the formula $\| (x, y) \|_p = (|x|^p + |y|^p)^{1/p}$. For $p = \infty$ it is defined as $\| (x, y) \|_\infty = \sup(|x|, |y|)$.
\end{dfn}

Let us note that $\| \cdotp \|_1 \geq \| \cdotp \|_2 \geq \| \cdotp \|_\infty$ and that the case $p = 2$ is isotropic. We aim to prove that for nonsmooth anisotropy the solutions need not be unique (and in general are not unique); to achieve this goal, we will study how do minimal surfaces with respect to the $p-$th norm look like. At first let us see an example that the solution is unique:

\begin{stw}\label{stw:uniqueness}
Let $\Omega \subset \mathbb{R}^2$ be an open, bounded, strictly convex set. Take $\phi(x, Du) = \| Du \|_1$. Let $f \in C(\partial \Omega)$. Denote by $u$ the solution to isotropic LGP for $f$. Then, if the boundaries of superlevel sets of $u$ are parallel to the axes of the coordinate system, then $u$ is a unique solution of the anisotropic LGP with respect to the $l^1$ norm.
\end{stw}

\begin{dd}
Let $v \in BV(\Omega)$, $Tv = f$. Then

\begin{equation*}
\int_\Omega |Dv|_1 \geq \int_\Omega |Dv|_2 \geq \int_\Omega |Du|_2.
\end{equation*}
By uniqueness of solution to Euclidean LGP the second inequality is strict, if only $u \neq v$. As the boundaries of superlevel sets of $u$ are parallel to the axes of the coordinate system, we have $\int_\Omega |Du|_1 = \int_\Omega |Du|_2$; it follows that $u$ is a unique solution to the anisotropic LGP. \qed
\end{dd}

\begin{prz}\label{ex:uniqueness}
Let $\Omega = B(0,1)$. Take $\phi(x, Du) = \| Du \|_1$. Let $f(\theta) = \cos(2 \theta)$. We construct the isotropic solution $u$ using Sternberg-Williams-Ziemer construction. We notice, as the picture below shows, that the boundaries of superlevel sets of $u$ are parallel to the axes of the coordinate system.

\begin{figure}[h]
\includegraphics[scale = 0.15]{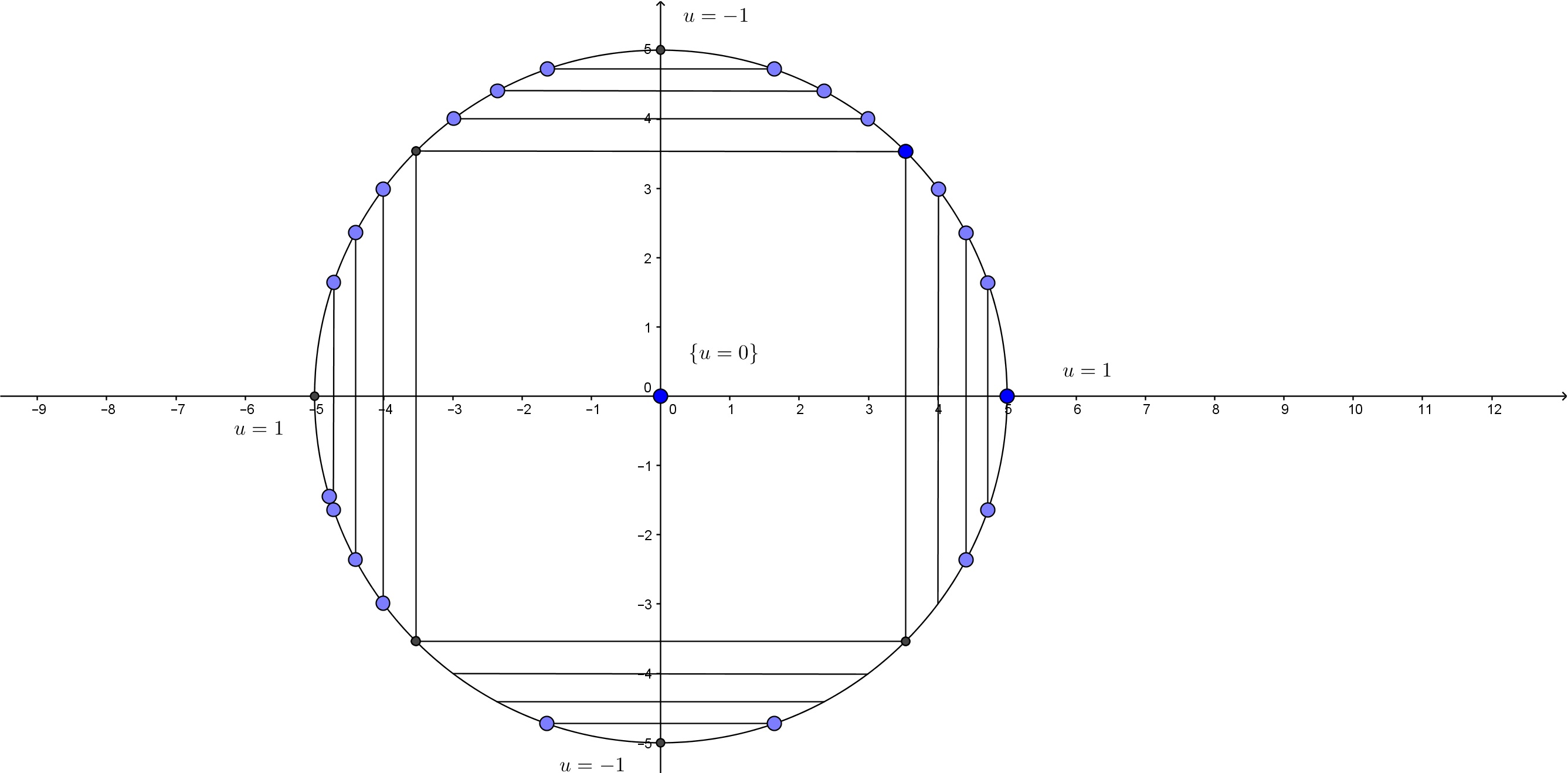}
\end{figure}
By Proposition \ref{stw:uniqueness} the solution to the anisotropic LGP is unique. \qed
\end{prz}

\begin{stw}\label{stw:niejednoznacznosc}
Let $\Omega \subset \mathbb{R}^2$ be an open, bounded, strictly convex set. Take $\phi(x, Du) = \| Du \|_1$. Let $f \in C(\partial \Omega)$. Denote by $u$ the solution to isotropic LGP for $f$. Then, if for some $t$ the boundaries of superlevel sets of $u$ are not parallel to the axes of the coordinate system, then the solution to the anisotropic LGP with respect to the $l^1$ norm is not unique.
\end{stw}

\begin{dd}
1. Take $v \in C^1(\Omega)$ with trace $f$. Then the co-area formula reads

\begin{equation*}
\int_\Omega |Dv|_1 = \int_{\mathbb{R}} P_1 (E_t, \Omega) dt,
\end{equation*}
in particular $v$ is a solution to anisotropic LGP iff $P_1(E_t, \Omega)$ is minimal for a.e. $t$. As $v$ is smooth, $v |_{\partial\Omega} = f$, then by Sard theorem for a.e. $t$ the set $\{ v = t \}$ is a smooth manifold; as such, it is an at most countable sum of smooth curves disjoint in $\Omega$. 

2. We want to find the lower bound for $\int_\Omega |Dv|_1$. We shall find it for a larger class of functions: continuous functions, for which the sets $\{ v = t \}$ are at most countable sums of smooth curves disjoint in $\Omega$. We have to extend our class of functions, as we need to be able to eliminate closed curves from the disjoint sum: if there were any closed curves, then by setting $v = t$ in the open set enclosed by such curves we obtain a function with strictly smaller total variation, but not necessarily smooth. Thus we may assume that $\partial \{ v \geq t \}$ is a disjoint sum of open curves. Let us note that they must end in points $p \in f^{-1}(t) \subset \partial\Omega$. 

3. According to the co-area formula, it is sufficient to construct superlevel sets of $v$ such that $P_1(E_t, \Omega)$ is minimal; then $\int_\Omega |Dv|_1$ would be minimal as well. Let us suppose additionally that $\partial E_t$ does not contain any vertical intervals, i.e. we may represent a level set from the point $(x,y)$ to $(z,t)$ as a graph of a $C^1$ function $g$. Let us note that at the point $((s, g(s)))$ the Radon-Nikodym derivative $\nu^{\chi_{E_t}}$ is perpendicular to the level set, so it is a vector $(- \sin \theta, \cos \theta)$, where $g'(s) = \tan \theta$. Thus $\phi(x, \nu^{\chi_{E_t}}) = |\sin \theta| + |\cos \theta|$. As $|D \chi_{E_t}| = \mathcal{H}^{n-1}|_{\partial E_t}$, then, using the representation introduced by Proposition \ref{stw:repcalkowa}, we have to minimize the integral (we may assume that $x < z$):

\begin{equation*}
P(E_t, \Omega) = \int_\Omega \phi(x, \nu^{\chi_{E_t}}) |D \chi_{E_t}| = \int_{\partial E_t} (|\sin \theta| + |\cos \theta|) d\mathcal{H}^{n-1} =
\end{equation*}
\begin{equation*}
= \int_{x}^{z} (|\sin \theta| + |\cos \theta|) \sqrt{1 + (\tan \theta)^2} dx = \int_{x}^{z} (|\sin \theta| + |\cos \theta|) \frac{1}{|\cos \theta|} dx = 
\end{equation*} 
\begin{equation*}
= \int_{x}^{z} (1 + |\tan \theta|) dx = |z - x| + \int_{x}^{z} |g'| dx \geq |z - x| + |t - y|,
\end{equation*}
where the inequality becomes equality iff $g$ is monotone (remember we assumed it to be $C^1$). Thus there are multiple functions minimizing this integral.

4. Now we allow $\partial E_t$ to contain vertical intervals. The difference is purely technical, as we have to divide our integral into two parts. Let us suppose that the (orientated) length of $i-$th vertical interval equals $\lambda_i$, then we have

\begin{equation*}
\int_{\text{graph part of } \partial E_t} (|\sin \theta| + |\cos \theta|) dl + \int_{\text{vertical part of } \partial E_t} (1 + 0) dl = \int_{x}^{z} (1 + |g'|) dx + \sum_{i = 1}^{\infty} |\lambda_i| = 
\end{equation*}
\begin{equation*}
= \int_{x}^{z} |g'| dx + |z - x| + \sum_{i = 1}^{\infty} |\lambda_i| \geq |t - y - \sum_{i = 1}^{\infty} \lambda_i| + |z - x| + \sum_{i = 1}^{\infty} |\lambda_i| \geq |z - x| + |t - y|,
\end{equation*}
where the inequality becomes equality iff $g$ is monotone (remember we assumed it to be $C^1$) and all the vertical intervals are orientated in the same direction as $g'$. Thus there are multiple functions minimizing this integral. We have proved that in a class containing all smooth functions the problem of minimizing perimeter of a set $E_t$ doesn't have a unique solution.

5. Let us denote by $u$ the solution to the Euclidean LGP. Let us notice that intervals are graphs of monotone functions, so an interval mimimizes the above integral; thus, by co-area formula, the value of $\int_\Omega |Dv|_1$ is bounded from below by
\begin{equation*}
\int_\Omega |Dv|_1 = \int_{\mathbb{R}} P_1(E_t, \Omega) \geq \int_{\mathbb{R}} P_1(\{ u > t \}, \Omega),
\end{equation*}
so by Remark \ref{lem:scisaniz} such inequality holds for all $v \in BV_1(\Omega)$ such that $Tv = f$. In particular, the Euclidean solution is also a solution to the anisotropic LGP. But if we choose $v$ such that its level sets $\{ v = t \}$ be monotone for almost all $t$, then its total variation is exactly the same (it is possible due to the non-parallelism assumption). Thus the solution to this anisotropic LGP is not unique. \qed
\end{dd}

\begin{prz}\label{ex:l1}
Let $\Omega = B(0,1)$. Take $\phi(x, Du) = \| Du \|_1$. Let $f \in C^{\infty}(\partial\Omega)$ be given as $f = \cos(2 \theta - \pi/2)$. Then the solution to the anisotropic LGP is not unique. At first, let us see that the Euclidean solution is a rotation of the function $u$ from Example \ref{ex:uniqueness}, so we may apply the procedure from Proposition \ref{stw:niejednoznacznosc}. We observe that for fixed $t \in (0,1)$ its preimage contains points of the form $A_1 = (a,b)$, $A_2 = (b, a)$, $A_3 = (-a, -b)$, $A_4 = (-b, -a)$; then, applying the above calculation to the function $f$, we see that the two possible connections, $A_1 A_2, A_3 A_4$ and $A_1 A_4, A_2 A_3$ have perimeter lengths $4 |a - b|$ and $4 |a + b|$ respectively; we choose the former as the level set $E_t$. Similar calculation holds for $t \in (-1,0)$. But if we choose $v$ such that its level sets $\{ v = t \}$ be monotone for almost all $t$, then their perimeter (and, by co-area formula, its total variation) stays exactly the same. Thus the solution to this anisotropic LGP is not unique; an example of a non-Euclidean solution is presented on the picture below.

\begin{figure}[h]
\includegraphics[scale = 0.15]{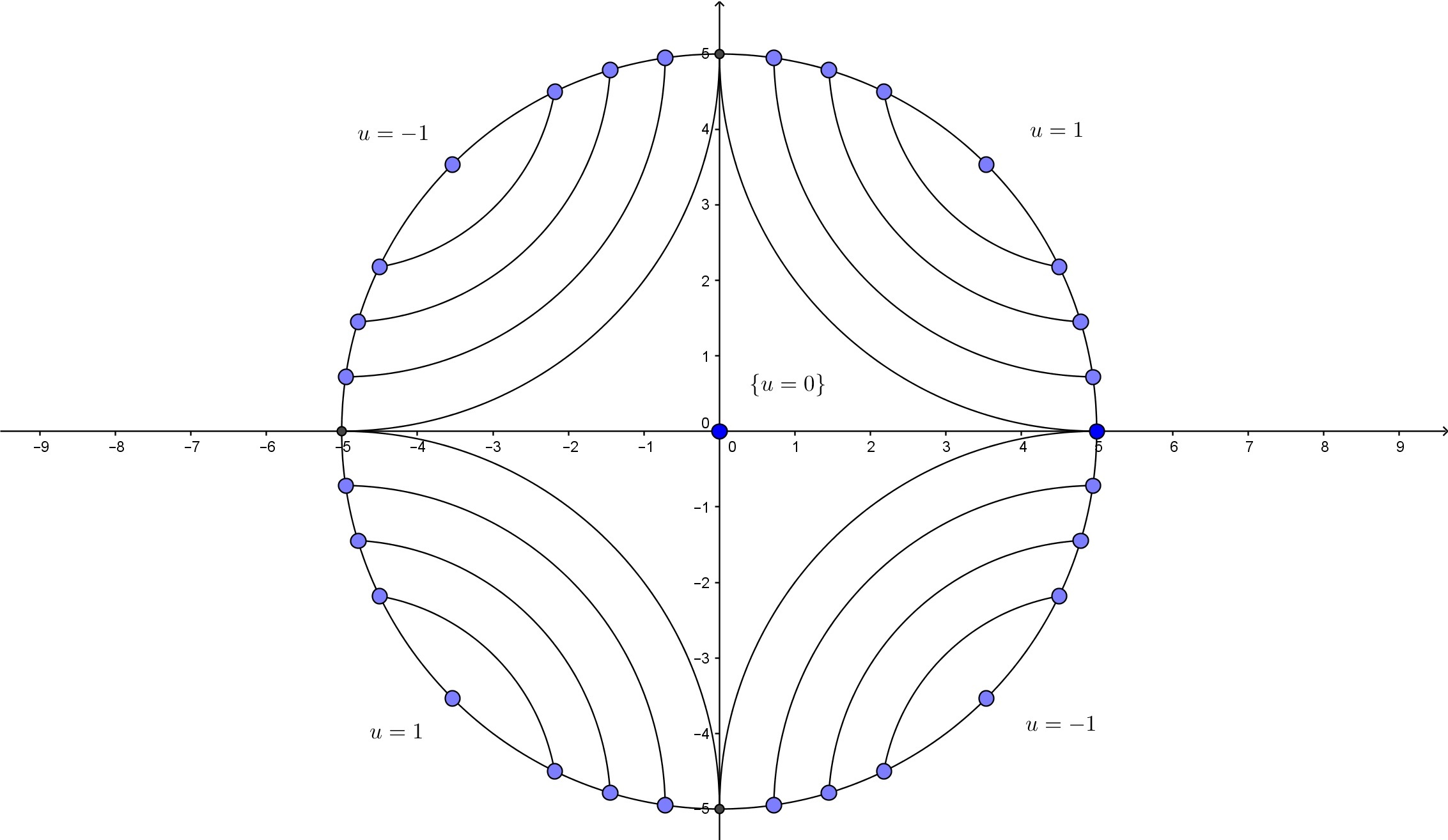}
\end{figure} \qed
\end{prz}

\begin{prz}\label{ex:linfty}
Now let $p = \infty$. If we make a similar calculation, we obtain that the perimeter of a level set connecting points $(x,y)$ with $(z,t)$ equals

\begin{equation*}
\int_{\partial E_t} \max(|\sin \theta|, |\cos \theta|) d\mathcal{H}^{n-1} = \int_{x}^{z} \max(|\sin \theta|, |\cos \theta|) \sqrt{1 + (\tan \theta)^2} dx = 
\end{equation*}
\begin{equation*}
= \int_{x}^{z} \max(|\sin \theta|, |\cos \theta|) \frac{1}{|\cos \theta|} dx = \int_{x}^{z} \max(1, |\tan \theta|) dx = \int_{x}^{z} \max(1, |g'|) dx \geq |z - x|,
\end{equation*}
where the inequality becomes equality iff $|g'| \leq 1$; in other words, the angle between the level set and the $x$ coordinate axis is not greater than $\frac{\pi}{4}$. Thus, if we take the function $f(\theta) = \cos(2 \theta)$, the solution is not unique; we apply this result for $t \in (-1,0)$ and then apply it again for $t \in (0,1)$ considering the level set as a function of $y$. A solution different than the Euclidean one is presented on the picture below. Nevertheless, it may still happen that the solution is unique: it is the case if we take such $f$ that the Euclidean solution has all level sets at an angle $\frac{\pi}{4}$ to the coordinate axes. For example we can take $f(\theta) = \cos(2 \theta - \frac{\pi}{2})$.

\begin{figure}[h]
\includegraphics[scale = 0.15]{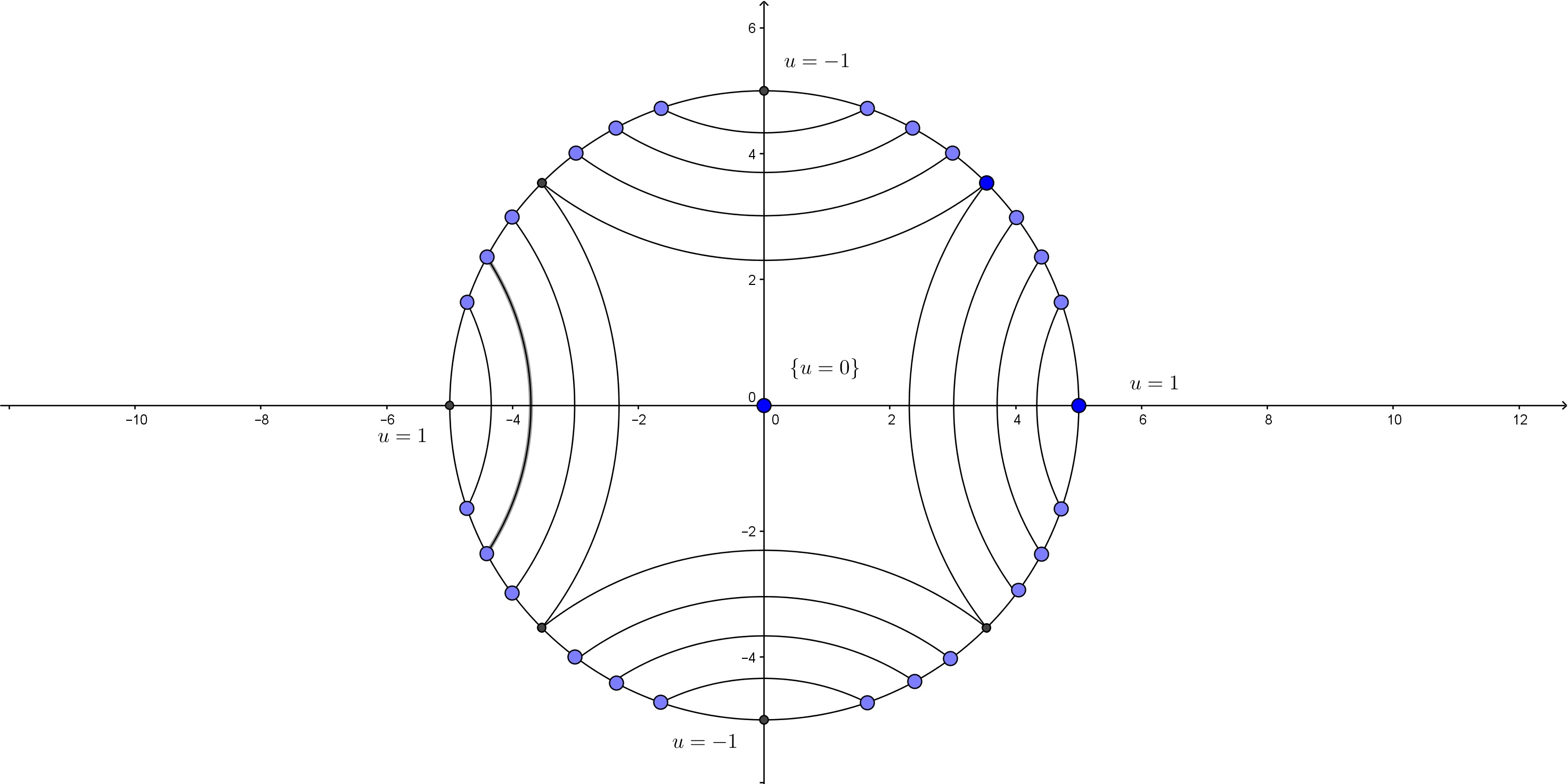}
\end{figure} \qed
\end{prz}

Now let $1 < p < \infty$. By \cite[Theorems 1.1, 1.2]{JMN} for continuous boundary data the anisotropic LGP has a unique solution, because the norm $\| \cdotp \|_p$ is a smooth function of the Euclidean norm outside $(0,0)$. We will show that connected components of boundaries of superlevel sets of functions of $\phi-$least gradient are line segments, similarly to for the isotropic norm $\| \cdotp \|_2$; in fact, due to an anisotropic analogue of Theorem \ref{twierdzeniezbgg} proved in \cite[Theorem 3.19]{Maz}, it is enough to show that the boundaries of minimal sets are line segments.

\begin{tw}\label{tw:anizotropia}
Let $\Omega \subset \mathbb{R}^2$ be an open convex set. Let the anisotropy be given by the function $\phi(x,Du) = \| Du \|_p$, where $1 < p < \infty$. Let $E$ be a $\phi-$minimal set with respect to $\Omega$, i.e. $\chi_E$ is a function of $\phi-$least gradient in $\Omega$. Then every connected component of $\partial E$ is a line segment.
\end{tw}

\begin{dd}
Let $(x,y),(z,t)$ be two points on the same connected component of $\partial E$. We have to minimize an integral analogous to the previous one (notation stays the same):

\begin{equation*}
L(x, g, g') = \int_{\partial E_t} (|\sin \theta|^p + |\cos \theta|^p)^{\frac{1}{p}} d\mathcal{H}^{n-1} = \int_{x}^{z} (|\sin \theta|^p + |\cos \theta|^p)^{\frac{1}{p}} \sqrt{1 + (\tan \theta)^2} dx = 
\end{equation*}
\begin{equation*}
= \int_{x}^{z} (|\sin \theta|^p + |\cos \theta|^p)^{\frac{1}{p}} \frac{1}{|\cos \theta|} dx = \int_{x}^{z} (1 + |\tan \theta|^p)^{\frac{1}{p}} dx = \int_{x}^{z} (1 + |g'|^p)^{\frac{1}{p}} dx.
\end{equation*}
The Euler$-$Lagrange equation for the functional $L$ takes form

\begin{equation*}
0 = \frac{\partial L}{\partial g} = \frac{d}{dx} (\frac{\partial L}{\partial g'}) = \frac{d}{dx} (\mathrm{sgn}(g') (g')^{p-1} (1 + |g'|^p)^{\frac{1}{p} - 1})
\end{equation*}
\begin{equation*}
\mathrm{sgn}(g') (g')^{p-1} (1 + |g'|^p)^{\frac{1}{p} - 1} = \text{ const}.
\end{equation*}
Taking absolute value and raising both sides to power $\frac{p}{p-1}$ we obtain

\begin{equation*}
\frac{|g'|^p}{1 + |g'|^p} = \text{ const} = C,
\end{equation*}
thus $g' =$ const. Thus the anisotropic minimal surface connecting points $(x,y)$ and $(z,t)$ is a line segment. \qed
\end{dd}

\begin{ak}
This paper is based on my master's thesis. My supervisor was Piotr Rybka, whom I would like to thank for many fruitful discussions on this paper. The author receives the WCNM scholarship. 
\end{ak}

\end{document}